\documentclass[amsmath,english,a4paper,graphicx,12pt]{article}
\usepackage{doc}
\usepackage{graphicx}
\usepackage{amsmath}
\usepackage{amsfonts}
\usepackage{amssymb}
\usepackage{babel}
\usepackage{latexsym}

\title{\bf On the Paley-Wiener theorem in the Mellin transform setting}
\author{Carlo Bardaro, \thanks{
Department of Mathematics and Computer Sciences, University of Perugia,
via Vanvitelli 1, I-06123 Perugia, Italy, e-mail: 
carlo.bardaro@unipg.it} \and
Paul L. Butzer, \thanks{Lehrstuhl A fuer Mathematik, RWTH Aachen, Templergraben 55, Aachen, D-52056, Germany}\and
Ilaria Mantellini, \thanks{Department of Mathematics and Computer Sciences, University of Perugia,
via Vanvitelli 1, I-06123 Perugia, Italy, e-mail: 
mantell@dmi.unipg.it}\and Gerhard Schmeisser\thanks{Department Mathematik, FAU Erlangen-N\"{u}rnberg, Cauerstr. 11, 91058 Erlangen, Germany, email: G.Schmeisser@gmx.de}
}
\begin{document}
\maketitle
\noindent
{\bf Abstract.} In this paper we establish a version of the Paley-Wiener theorem of Fourier analysis in the frame of the Mellin transform. We provide two different proofs, one involving complex analysis arguments, namely the Riemann surface of the logarithm and Cauchy theorems, and the other one employing a Bernstein inequality here derived for Mellin derivatives.

\section{\bf Introduction}
\vskip0.5cm
An entire function $f: \mathbb{C} \rightarrow \mathbb{C}$, bounded on the real line, is said to be of exponential type $T >0$, if 
it satisfies the inequality
\begin{eqnarray}
|f(z)| \leq C \exp(T |\Im z|)\quad\quad\quad (z \in \mathbb{C}),
\end{eqnarray}
with $\Im z$ being the imaginary part of the complex number $z.$

The Bernstein space $B^2_T$ consists of all functions $f$ of exponential type $T,$ whose restriction to $\mathbb{R}$ belongs 
to $L^2(\mathbb{R}).$ The essentially deep assertion of the Paley-Wiener theorem states that if $f \in B^2_T$ for some $T >0$ then, denoting as usual by 
$\widehat{f}$ the $L^2-$Fourier transform of $f,$ one has $\widehat{f}(v) = 0$ a.e. outside 
the interval $(-T, T),$ thus
\begin{eqnarray}
f(z) = \frac{1}{\sqrt{2 \pi}}\int_{-T}^T \widehat{f}(v) e^{ivz}dv \quad \quad \quad (z \in \mathbb{C}).
\end{eqnarray}
Conversely, employing elementary arguments of complex analysis, every function $f \in L^2(\mathbb{R})$ such that $\widehat{f}(v) = 0$ a.e. outside 
the interval $(-T, T),$ has an extension to the complex plane given by (2), which belongs to $B^2_T.$

The Paley-Wiener theorem is a basic result for many results in Fourier analysis, in particular the Shannon sampling theorem for functions belonging to the Bernstein space 
(see \cite{BSS} and references therein).

The above well-known result and its diverse distributional versions can be found in \cite{RU}, \cite{RU1}, \cite{YO}. Their proofs  usually employ complex analysis arguments, like Cauchy and Morera theorems.

The aims of this paper is to give a direct extension of the Paley-Wiener theorem in Mellin transform setting using two different approaches. An important application of this extension is the mathematical foundation of the exponential sampling theorem of signal analysis. 

For a given $c \in \mathbb{R},$ let $X^2_c$ be the space of all functions $f: \mathbb{R}^+ \rightarrow \mathbb{C},$ 
such that $f(\cdot) (\cdot)^{c -1/2} \in L^2(\mathbb{R}^+).$

Denoting by $M^2[f]\equiv [f]^\wedge_{M^2},$ the Mellin transform of a function $f \in X^2_c$ (see Section 2), we say that $f$ is Mellin bandlimited to an interval $[-T,T]$ if  
$[f]^\wedge_{M^2}(c+it) = 0$ whenever $|t| > T.$ 

We show that a Mellin bandlimited function $f$ cannot be extended to the whole complex plane as an entire function, but it can be extended to the Riemann surface of the logarithm as an analytic function. Using this fact we define the counterparts of the Bernstein spaces of Fourier analysis in the Mellin transform setting (the Mellin-Bernstein spaces) and we prove an extension of the Paley-Wiener theorem, employing complex variable arguments (Section 4). 

In Section 5 we state a Mellin version of the classical Bernstein inequality (see e.g., \cite{BO}, \cite{RS}), involving Mellin 
derivatives. Such a result was already given in \cite{BJ3} for Mellin bandlimited functions. However the present proof is fully independent of the Paley-Wiener theorem and it uses only the properties of the Mellin Bernstein spaces. Also, it provides a sharp inequality. Employing this inequality, we give an alternative proof of the Paley-Wiener theorem in Mellin transform setting, which avoids the use of the Riemann surfaces. This approach is based on a technique developed in \cite{AN}, which characterizes the support of a function $f.$ 
In Section 6 we apply the result to the exponential sampling, which arises in optical physics and engineering (see e.g. \cite{BP}).
 
 Mathematical versions of the exponential sampling theorem
(see Section 6), were given in \cite{BJ3}, \cite{BJ4}, \cite{BBM1}. Our Paley-Wiener theorem enables one to obtain the esponential sampling formula directly for functions in the Mellin-Bernstein spaces. Also, the same holds for another very connected result, namely the Mellin reproducing kernel formula, established in \cite{BBM1} for the space $B^2_T$ and not for $\widetilde{B}^2_T$ as below (see Definitions 1 and 2).

More recently, some different "real" versions of the Paley-Wiener theorem were stated, involving other approaches, (see e.g. \cite{AN}, \cite{AD}). These ``real'' approaches were also used in an earlier paper by V.K. Tuan in \cite{TU} for multidimensional Mellin (or inverse Mellin) transforms. 
We recall that Mellin transform analysis was recently widely developed in connection with various fields of mathematical analysis,  in particular in fractional calculus (see e.g. \cite{BKT},  \cite{BBM2}) and quadrature formulas on semi-infinite intervals (see e.g. \cite{SCH}).

\section{\bf Basic notions and preliminary results}

Let $C(\mathbb{R}^+)$ be the space of all continuous functions defined on $\mathbb{R}^+,$ and $C^{(r)}(\mathbb{R}^+)$ be the space of all functions in $C(\mathbb{R}^+)$ with $r-$order derivative in $C(\mathbb{R}^+).$ Analogously, by $C^\infty(\mathbb{R}^+)$ we denote the space of all infinitely differentiable functions.
By $L^1_{\mbox{loc}}(\mathbb{R}^+)$ we denote the space of all measurable functions which are integrable on every bounded interval in $\mathbb{R}^+.$

For $1\leq p < +\infty,$ let $L^p= L^p(\mathbb{R}^+)$~ be the space of all the Lebesgue measurable and $p-$integrable complex-valued functions defined on $\mathbb{R}^+$ endowed with the usual norm $\|f\|_p.$ Analogous notations hold for functions 
defined on $\mathbb{R}.$

For $p=1$ and $c \in \mathbb{R},$ let us consider the space
$$X_c = \{ f: \mathbb{R}^+\rightarrow \mathbb{C}: f(\cdot) (\cdot)^{c-1}\in L^1(\mathbb{R}^+) \}$$
endowed with the norm
$$ \| f\|_{X_c} = \|f(\cdot) (\cdot)^{c-1} \|_1 = \int_0^\infty |f(u)|u^{c-1} du.$$

More generally let $X^p_c$  denote the space of all  functions $f: \mathbb{R}^+\rightarrow \mathbb{C}$ such that $f(\cdot) (\cdot)^{c-1/p}\in L^p (\mathbb{R}^+),$ with $1<p< \infty.$ In an equivalent form, $X^p_c$ is the space of all functions $f$ such that $(\cdot)^c f(\cdot) \in L^p_\mu(\mathbb{R}^+),$ where $L^p_\mu= L^p_\mu(\mathbb{R}^+)$ denotes the Lebesgue space with respect to the (invariant) measure $\mu (A) = \int_A dt/t,$ for  any measurable set $A \subset \mathbb{R}^+.$

 The Mellin translation operator $\tau_h^c$, for $h \in \mathbb{R}^+,~c \in \mathbb{R},$~$f: \mathbb{R}^+ \rightarrow \mathbb{C},$ is denoted by
$$(\tau_h^c f)(x) := h^c f(hx)~~(x\in \mathbb{R}^+).$$
\noindent
Setting $\tau_h:= \tau^0_h,$ then  $(\tau_h^cf)(x) = h^c (\tau_hf)(x),$  
$\|\tau_h^c f\|_{X_c} = \|f\|_{X_c}.$
\vskip0.3cm
In the Mellin frame, the natural concept of a  pointwise derivative of a function $f$ is given by the limit of the difference quotient involving the Mellin translation; thus if $f'$ exists,
$$\lim_{h \rightarrow 1}\frac{\tau_h^cf(x) - f(x)}{h-1} =  x f'(x) + cf(x).$$
This gives the motivation for the following definition (see \cite{BJ}):
the pointwise Mellin differential operator $\Theta_c,$ or the pointwise Mellin derivative $\Theta_cf$ of a function $f: \mathbb{R}^+ \rightarrow \mathbb{C}$ and $c \in \mathbb{R},$ is defined by
\begin{eqnarray}
\Theta_cf(x) := x f'(x) + c f(x),~~x \in \mathbb{R}^+
\end{eqnarray}
provided $f'$ exists a.e. on $\mathbb{R}^+.$ The Mellin differential operator of order $r \in \mathbb{N}$ is defined recursively by
\begin{eqnarray}
\Theta^1_c := \Theta_c ,\quad\quad \Theta^r_c := \Theta_c (\Theta_c^{r-1}).
\end{eqnarray}
For convenience, set $\Theta^r:= \Theta^r_0$ for $c=0$ and $\Theta_c^0 := I,$ $I$ denoting the identity operator.
For instance, the first three Mellin derivatives are given by:
\begin{eqnarray*}
\Theta_cf(x) &=& xf'(x) + cf(x),\\
\Theta^2_cf(x) &=& x^2 f''(x) + (2c+1) xf'(x) + c^2f(x),\\
\Theta^3_cf(x) &= &x^3 f'''(x) + (3c+3)x^2f''(x) \\&& +(3c^2 + 3c +1)xf'(x) + c^3 f(x).
\end{eqnarray*}
A definition of the Mellin derivative was given also in \cite{MA} in case $c=0$, in a slight different, but equivalent form, using the quotients
$$\frac{f(xh^{-1}) - f(x)}{\log h}.$$
\vskip0,4cm
In order to define the underlying spaces such that $\Theta f$ is meaningful, recall that the set $AC_{\rm{loc}}(\mathbb{R}^+)$ of all locally absolutely continuous functions on $\mathbb{R}^+$, can be characterized as the space of all functions $f:\mathbb{R}^+ \rightarrow \mathbb{C}$ for which there exists a locally integrable function $g \in L^1_{\rm{loc}}(\mathbb{R}^+)$ and a constant $\alpha \in \mathbb{C}$ such that 
$$f(x) = \alpha + \int_1^x g(u)du ~~~(x \in \mathbb{R}^+).$$
Thus for $f \in AC_{\rm{loc}}(\mathbb{R}^+)$ the derivative of $f$ exists a.e. with $f' = g$ a.e.
For any $r \in \mathbb{N},$ and $p\geq 1,$ we define the Mellin-Sobolev space $W^{p,r}_c$ as the space of all functions 
$f \in X^p_c$ such that there exists $g \in C^{(r-1)}(\mathbb{R}^+)$ such that $f=g$ a.e. and $g^{(r-1)}\in AC_{\rm{loc}}(\mathbb{R}^+)$ and   $ \Theta^r_cg\in X^p_c.$

The Mellin transform of a function $f\in X_c$ is defined by (see e.g. \cite{MA}, \cite{GPS},  \cite{BJ2})
$$ M[f](s) \equiv [f]^{\wedge}_M (s) = \int_0^\infty u^{s-1} f(u) du,~~(s=c+ it, t\in \mathbb{R}).$$
Basic properties of the Mellin transform are the following
$$M[af(\cdot) + bg(\cdot)](s) = a M[f](s) + bM[g](s)$$
for $f,g \in X_c,~a,b \in \mathbb{R}$ and
$$|M[f(\cdot)](s)| \leq \|f\|_{X_c}~~(s = c+it).$$
The inverse Mellin transform $M^{-1}_c[g]$ of the function $g \in L^1(\{c\} \times i \mathbb{R}),$ is defined by:
\begin{eqnarray*}
&&M^{-1}_c[g](x) \equiv M^{-1}_c[g(c+it)](x) :=\\
&& \frac{x^{-c}}{2 \pi}\int_{-\infty}^{+\infty} g(c+it) x^{-it}dt, ~~(x \in \mathbb{R}^+),
\end{eqnarray*}
where  $L^p(\{c\} \times i \mathbb{R}),$ for $p \geq 1,$ will mean the space of all functions $g:\{c\} \times i \mathbb{R} \rightarrow \mathbb{C}$ with 
$g(c +i\cdot) \in L^p(\mathbb{R}^+).$
\vskip0,3cm
We have the following preliminary results (see \cite{BJ2}, \cite{BBM1})
\vskip0,3cm
\newtheorem{Lemma}{Lemma}
\begin{Lemma}~[Inversion Theorem in $X_c$]. 
If $f \in X_c$ is such that $M[f] \in L^1(\{c\} \times i \mathbb{R}),$ then 
$$M_c^{-1}[M[f]](x) = \frac{x^{-c}}{2 \pi} \int_{-\infty}^\infty [f]^\wedge_M(c+it)x^{-it}dt = f(x) \quad \quad (a.e.~ x \in \mathbb{R}^+).$$
\end{Lemma}
\vskip0,3cm
\noindent
 Under the hypothesis that $f \in X_c$ and $M[f] \in L^1(\{c\} \times i \mathbb{R}),$ the following lemma will enable us to work  in a practical Hilbert space setting.
\begin{Lemma}
 If $f \in X_c$  and $M[f] \in L^1(\{c\} \times i \mathbb{R}),$ then 
$f \in X^2_c.$ 
\end{Lemma}
\vskip0,4cm

More generally, for $1<p \leq 2,$  the Mellin transform $M^p$ of $f \in X^p_c,$ is given by (see \cite{BJ4})
$$M^p[f](s) \equiv [f]^{\wedge}_{M^p} (s) = \mbox{l.i.m.}_{\rho \rightarrow +\infty}~\int_{1/\rho}^\rho f(u) u^{s-1}du,$$
for $s=c+it,$ 
in the sense that
$$\lim_{\rho \rightarrow \infty}\bigg\|M^p[f](c+it) - \int_{1/\rho}^\rho f(u) u^{s-1}du\bigg\|_{L^p(\{c\}\times i \mathbb{R})} = 0.$$
In the following we are interested to the case $p=2.$ 

Analogously, we define the inverse Mellin transform of a function $g \in X^2_c$ by 
$$M^{2, -1}_c[g](s) = \mbox{l.i.m.}_{\rho \rightarrow +\infty}~\frac{1}{2 \pi}\int_{1/\rho}^\rho g(c+it) x^{-c-it}dt,$$
and the inverse formula holds for any $f \in X^2_c$ (see \cite{BJ4})
$$M^{2, -1}_c[M^2[f]](x) = f(x),\quad \quad a.e.~ x \in \mathbb{R}^+.$$

For functions in $X_c \cap X^2_c,$ we have the following important "consistency" property of the Mellin transform (see \cite{BJ4}):
\vskip0,3cm
\begin{Lemma} 
If $f \in X_c \cap X^2_c$ then the Mellin transforms $M[f]$ and $M^2[f]$ coincide, i.e. $M[f](c+it) = M^2[f](c+it)$ for almost all $t \in \mathbb{R}.$
\end{Lemma}
\vskip0,4cm
Moreover, the following Mellin version of the Plancherel Theorem holds (see \cite{BJ4}, Lemma 2.6)
\vskip0,3cm
\begin{Lemma}
The operator $M^2$ from $X^2_c$ onto $L^2(\{c\}\times i \mathbb{R})$ is bounded and norm preserving i.e. for $f \in X^2_c$
$$\|f\|_{X^2_c}= \frac{1}{\sqrt{2\pi}}\|M^2[f]\|_{L^2(\{c\}\times i \mathbb{R})}.$$
\end{Lemma}
\vskip0,3cm
\begin{Lemma} If $f \in W^{2,r}_c,$ $r \in \mathbb{N},$ then 
\begin{eqnarray}
M^2[\Theta^r_cf](s) = (-it)^r M^2[f](s) \quad \quad (s = c +it, ~t \in \mathbb{R}).
\end{eqnarray}
\end{Lemma}
{\bf Proof}.
The proof can be found in \cite{BJ}, Proposition 6, for functions belonging to $W^{1,r}_c,$ but using formula (8.7) in \cite{BJ} and the convolution theorem for the Mellin transform $M^2,$ (see \cite{BJ4}, Lemma 2.7), it can be extended to the space 
$W^{2,r}_c.~~\Box$

\section{Mellin bandlimited functions and their properties}
We begin by introducing the space of the Mellin bandlimited functions (see \cite{BJ3}, \cite{BJ4})
\newtheorem{Definition}{Definition}
\begin{Definition}\label{def1}
Let $B^{2}_{c, T}$ denote the space of all functions in $X_c^2\cap C(\mathbb{R}^+)$ 
such that $[f]^\wedge_{M^2}(c+it) = 0$ a.e. for $|t| > T.$ Analogously, by $B^1_{c,T}$ we denote the space of all functions in $X_c^1\cap C(\mathbb{R}^+)$ 
such that $[f]^\wedge_{M}(c+it) = 0$ for all $|t| > T.$
\end{Definition}
By Lemma 2 we have $B^1_{c,T} \subset X^2_c$ and 
it is easily seen that $B^1_{c,T} \subset B^2_{c, T}.$ 

Using now the inverse Mellin transform $M^{2, -1}_c$, for any $f \in B^2_{c,T}$ we have the representation
$$f(x) = \frac{1}{2 \pi}\int_{-T}^T [f]^\wedge_{M^2}(c+it) x^{-c - it}dt,$$
or
\begin{eqnarray}
x^c f(x) = \frac{1}{2 \pi}\int_{-T}^T [f]^\wedge_{M^2}(c+it) e^{-it\log x}dt,\quad \quad x>0.
\end{eqnarray}

The following result states the behaviour of a Mellin bandlimited function at the origin and at infinity. 
\vskip0,2cm
\newtheorem{Theorem}{Theorem}
\begin{Theorem}
If $f \in B^2_{c,T}$ then $f \in C^\infty(\mathbb{R}^+)$ and for any $k=0,1,\ldots$
$$\lim_{x \rightarrow 0^+} x^c\Theta^k_cf(x) = \lim_{x \rightarrow +\infty} x^c\Theta^k_cf(x) = 0.$$
\end{Theorem}

{\bf Proof}. Putting $F(t) := [f]^\wedge_{M^2}(c+it),$ we obviously have $F \in L^1(\mathbb{R})$ and therefore using the inversion Mellin 
transform, we can write
$$f(x) = \frac{1}{2\pi}\int_{-T}^T F(t) x^{-c-it}dt = \frac{x^{-c}}{2 \pi} \int_{-T}^T F(t) e^{-it \log x}dt,$$
or
\begin{eqnarray}
x^c f(x) =  \frac{1}{2 \pi} \int_{-T}^T F(t) e^{-it \log x}dt.
\end{eqnarray}
Hence, putting $g(x):= x^c f(x),$ $x>0,$ by the classical Riemann-Lebesgue lemma (see e.g. \cite{BN}, page 51, Exercise 8(i)) we get
$$\lim_{x \rightarrow 0^+} x^c f(x) = \lim_{x \rightarrow +\infty} x^c f(x) = 0.$$
Now, under differentiation, we have also
$$g'(x) = \frac{-i}{2 \pi}\int_{-T}^T tF(t) e^{-it \log x}\frac{dt}{x},$$
or
$$xg'(x) = \frac{-i}{2 \pi}\int_{-T}^T tF(t) e^{-it \log x} dt,$$
and so, since also the function $t F(t)$ is integrable over $[-T,T],$ using again the Riemann-Lebesgue lemma,
we obtain
$$\lim_{x \rightarrow 0^+}\Theta g(x) = \lim_{x \rightarrow +\infty} \Theta g(x) = 0.$$
Since $x^c \Theta_c f(x)= \Theta g(x)$ we obtain the assertion for $k=1.$ 
The general case follows by the same method and the representation (see \cite{BBM2})
$$x^c \Theta_c^kf(x) = \Theta^k g(x).~\Box$$
\vskip0,4cm
In particular for $c=0,$ we have $g=f$ and for $k=0,1,\ldots,$
$$\lim_{x \rightarrow 0^+} \Theta^k f(x) = \lim_{x \rightarrow + \infty} \Theta^k f(x) = 0,$$
and as a consequence
$$\lim_{x \rightarrow 0^+}x^k f(x) = \lim_{x \rightarrow + \infty}x^k f(x) = 0.$$
\vskip0,4cm
\noindent
We now extend the definition of the integral in (6) to $z \in \Omega,$ with 
$\Omega:=\{z \in \mathbb{C}: z \not\in \mathbb{R}^+_0\},$ i.e. 
$$z^c f(z) = \frac{1}{2 \pi}\int_{-T}^T [f]^\wedge_{M^2}(c+it) e^{-it\log z}dt, \quad, \quad z \in \Omega,$$
which is an analytic function on $\Omega$ (this may be proved using the classical Morera theorem, as in the classical instance). 
Now, $\log z = \log |z| + i \arg (z),$ where $0 \leq \arg (z) < 2\pi.$ Therefore $e^{-it \log z} = e^{-it \log |z|}e^{t \arg (z)}.$ Thus:
\begin{eqnarray*}
&&|z^c f(z)| \leq \frac{1}{2 \pi}\int_{-T}^T |[f]^\wedge_{M^2}(c+it)| e^{t \arg (z)}dt \\ 
&&\leq \frac{e^{T |\arg (z)|}}{2 \pi}\int_{-T}^T |[f]^\wedge_{M^2}(c+it)|dt = C e^{T |\Im (\log z)|}.
\end{eqnarray*}

We now prove that the function $g(z) := z^c f(z)$ cannot be extended to the whole complex plane as an entire function.
\vskip0,3cm
\begin{Theorem}
Let $f \in B^2_{c,T}$ be a non zero function. Then $g $ cannot be extended to the whole complex plane as an entire function.
\end{Theorem}
{\bf Proof}. We have shown before that 
\begin{eqnarray}
z^c f(z) = \frac{1}{2 \pi}\int_{-T}^T [f]^\wedge_{M^2}(c+it) e^{-it\log z}dt, \quad \quad z \in \Omega,
\end{eqnarray}
is analytic over $\Omega.$ Let now $x \in ]0, \infty[$ be fixed. We define
$$g(x_{\pm}) := \lim_{\varepsilon \rightarrow 0^+}
\frac{1}{2 \pi}\int_{-T}^T [f]^\wedge_{M^2}(c+it) e^{-it\log (x\pm i\varepsilon)}dt.$$
We have,
\begin{eqnarray*}
\lim_{\varepsilon \rightarrow 0^+}e^{-it\log (x\pm i\varepsilon)} =
\left\{\begin{array}{ll} e^{-it\log x}~~&\mbox{for the upper sign}\\ \\
e^{-it\log x}e^{2\pi t}~~&\mbox{for the lower sign}.
\end{array}
\right.
\end{eqnarray*}
Thus, interchanging the limit and the integration, if the right hand side of (8) has an analytic continuation to the whole $\mathbb{C},$ we must have $g(x_+) = g(x_-),$ or equivalently 
$$\frac{1}{2 \pi}\int_{-T}^T [f]^\wedge_{M^2}(c+it) e^{-it\log x} e^{\pi t}
\sinh (\pi t) dt = 0.$$
If $x$ runs through $\mathbb{R}^{+},$ then  $u:= \log x$ runs through the whole of $\mathbb{R}.$ Since $f$ is Mellin bandlimited to $[-T,T],$ one may extend the range of integration to $\mathbb{R},$, i.e.
$$\frac{1}{2 \pi}\int_{-\infty}^\infty [f]^\wedge_{M^2}(c+it) e^{-itu} 
e^{\pi t}\sinh (\pi t) dt = 0,$$
for all $u \in \mathbb{R}.$ This means that the Fourier transform of $[f]^\wedge_{M^2}(c+it) e^{\pi t} \sinh (\pi t)$ is identically zero. We conclude that $[f]^\wedge_{M^2}(c+it) $ must be zero almost everywhere, and so $f$ is identically zero by Lemma 4 and the continuity of $f.$ This is a contradiction. $\Box$
\vskip0,3cm
\noindent
From the above result one can deduce also, that a (classical) entire function cannot be Mellin-bandlimited unless it is identically zero. 

\section{The Paley-Wiener theorem for the Mellin transform}

We have seen that a Mellin bandlimited function cannot be extended to the whole complex plane as an entire function. However, as 
we shall prove, it has an analytic extension
on the Riemann surface of the logarithm. 
For establishing a Paley-Wiener theorem,
we first introduce the following Mellin-Bernstein space $\widetilde{B}^2_{c,T}$.
\begin{Definition}\label{def2}
The Mellin-Bernstein space $\widetilde{B}^2_{c,T}$ comprises all functions $f\in X_c^2$
for which $g(x):=x^cf(x)$ has an analytic extension on the Riemann
surface of the logarithm. Then, for each $k\in\mathbb{Z}$ this extension has
an analytic branch $g_k$ on $\Omega$ such that the following holds:
\begin{enumerate}
\item[(i)]
For $x>0$ the limits
$$ g_k^+(x):= \lim_{\varepsilon \rightarrow 0+} g_k(x+i\varepsilon)
\quad \hbox{and}\quad
g_k^-(x):= \lim_{\varepsilon \rightarrow 0+} g_k(x-i\varepsilon)$$
 exist and
\begin{eqnarray}
g_k^-(x)= g_{k+1}^+(x), \quad g_0^+(x)= g(x)
\quad (x>0).
\end{eqnarray}
\item[(ii)]
For $x>0$, let $U_x$ be an open disk in the right half-plane
with center at $x$. Then $\psi_k\,:\, U_x\rightarrow \mathbb{C}$ with
\begin{eqnarray*}
 \psi_k(z):= \left\{
\begin{array}{lll}
g_k(z) & \hbox{  for } & z\in U_x, \, \Im z <0,\\
g_k^-(z) & \hbox{ for } & z\in U_x \cap\mathbb{R},\\
g_{k+1}(z) & \hbox{ for } & z\in U_x, \,  \Im z>0
\end{array}
\right.
\end{eqnarray*}
is analytic.
\end{enumerate}
In addition, we require that the branches $g_k$ have the following
properties:
\begin{enumerate}
\item[(iii)]
There exists a constant $C>0$ such that for all $k\in\mathbb{Z}$ and
$\theta\in [0,2\pi]$
\begin{eqnarray}
|g_k(re^{i\theta})|\,\leq \, C e^{T|2\pi k +\theta|} \qquad (r>0).
\end{eqnarray}
\item[(iv)]
For $\theta\in [0,2\pi]$, we have
\begin{eqnarray}
\lim_{r\rightarrow 0} g_k(re^{i\theta})\,=\, \lim_{r\rightarrow \infty}
g_k(re^{i\theta})\,=\,0
\end{eqnarray}
uniformly with respect to $\theta$.
\end{enumerate}
In (10) and (11) the value $g_k(re^{i\theta})$ has to be
defined as $g_k^+(r)$ when $\theta=0$ and as $g_k^-(r)$ when
$\theta=2\pi$.
\end{Definition}

Now the desired Paley-Wiener theorem for the Mellin transform can be
stated as follows:
\begin{Theorem}[Paley-Wiener] $\quad \widetilde{B}^2_{c,T} \,=\,
B^2_{c,T}$.
\end{Theorem}
\noindent
{\bf Proof.\,}
Let $f$ be any function from $B^2_{c,T}$. Then, by the inversion formula (6), we have
$$ h(x):= e^{cx} f(e^x)\,=\, 
\frac{1}{2\pi} \int_{-T}^T [f]^\wedge_{M^2}(c+it) e^{-itx} \,dt \qquad (x>0).$$
Now it follows from the classical Paley-Wiener theorem that $h$ has an
extension to an entire function such that
\begin{eqnarray}
|h(x+iy)| \leq C e^{T|y|}\qquad (x, y\in\mathbb{R})
\end{eqnarray}
with some constant $C>0$ and  
\begin{eqnarray}
\lim_{x\rightarrow \pm\infty} h(x+iy)\,=\,0
\end{eqnarray}
uniformly with respect to $y$ on bounded intervals.

Formally, we have
$$ g(z)\,:=\, z^c f(z)\,=\, h(\log z).$$
This relation certainly holds for $z\in \mathbb{R}^+$. The logarithm can be
extended to its Riemann surface $S_{\log}$, and since $h$ is an entire 
function, this entails an extension of $g$ on $S_{\log}$.
On $\Omega$ the extension of the logarithm has analytic branches
$$\log_k(re^{i\theta})\,:=\, \log r + i(2\pi k+\theta)
\qquad (r>0, \,\, \theta\in (0, 2\pi), \,\, k\in\mathbb{Z}).$$
They induce branches $g_k$ of $g$, defined by
$$g_k(re^{i\theta})\,:=\, h(\log_k(re^{i\theta}))\,=\,
h(\log r + i(2\pi k+\theta)).$$
Now it is easily verified that statements (i) and (ii) of Definition~\ref{def2}
hold. Statements (iii) and (iv) are consequences of (12) and
(13), respectively. Altogether, we have shown that
$f\in \widetilde{B}^2_{c,T}$ and so $B^2_{c,T} \subset
\widetilde{B}^2_{c,T}.$

Now we prove the converse inclusion.
 We have to show that $[f]^\wedge_{M^2}(c+it) = 0$ for $|t| > T.$  For any $k \in \mathbb{Z}$ and $ z = re^{i \theta}$ with 
$0 \leq \theta < 2 \pi,$ we define
$$F_k(z) : = \frac{g_k(z)}{z}\exp(it(\log |z| + i(2 \pi k - \theta))),$$
which is analytic in $\Omega$ since the argument of the exponential function is an analytic branch of $it \log z.$ Hence the integral of $F_k$ along 
the closed curve $\gamma(R)$ as shown in Fig. 1 vanishes. 
 \begin{figure}[h!]
\begin{center}
\includegraphics[scale=0.4]{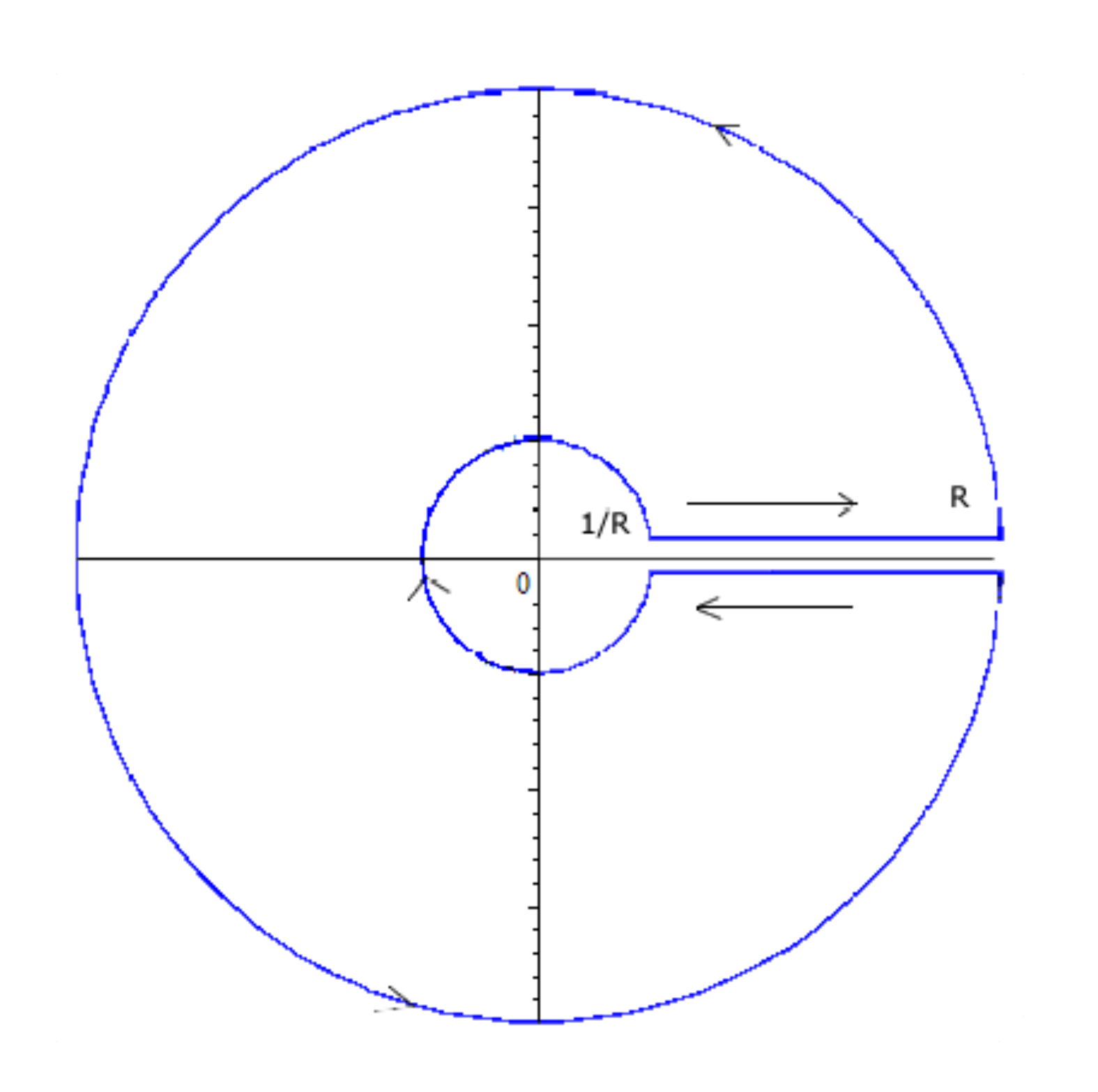}
\caption{The closed curve $\gamma(R)$}\label{fig:1}
\end{center}
\end{figure}

We may even let the two line segments of $\gamma(R)$ coincide on $\mathbb{R}$ provided we replace 
the integrand by
$$F_k^+(x) : = \frac{g_k^+(x)}{z}\exp(it(\log x + i2k\pi)),$$
for the integration from $1/R$ to $R$ and  by 
\begin{eqnarray*}
F_k^-(x) &:= &\frac{g_k^-(x)}{z}\exp(it(\log x + i2(k+1)\pi)) \\&= &
 \frac{g_{k+1}^+(x)}{z}\exp(it(\log x + i2(k+1)\pi)) = F^+_{k+1}(x),\end{eqnarray*}
 for the integration  from $R$ to $1/R.$
 Now for any $\rho >0$ let us consider a circle $\lambda (\rho)$ given by $z = \rho e^{i\theta},$ for $0 \leq \theta \leq 2 \pi,$ and let us define
 $$I_k(\rho):= \int_{\lambda(\rho)}F_k(z)dz = i \int_0^{2 \pi} g_k(\rho e^{i\theta}) e^{it \log \rho} e^{-t(2k\pi + \theta)}d\theta$$
 and 
 $$J_k(\rho) := \int_{1/\rho}^\rho F_k^+(u) du = \int_{1/\rho}^{\rho} g_k^+(u) e^{it \log u} e^{-t 2k\pi}\frac{du}{u}.$$
 Now we have
 $$0 = \int_{\gamma(R)} F_k(z) dz = J_k(R) + I_k(R) - J_{k+1}(R) - I_k(1/R)$$
 and so 
 \begin{eqnarray}
 J_k(R) - J_{k+1}(R) = I_{k}(1/R) - I_k(R).
 \end{eqnarray}
 Summing over $k = 0,1,\ldots n-1$ we get by (9)
 \begin{eqnarray}
 \int_{1/R}^R g(u) e^{it \log u}\frac{du}{u} = J_n(R) + \sum_{k=0}^{n-1}(I_k(1/R) - I_k(R)).
 \end{eqnarray}
 First let $t >T$ and let $\epsilon >0$ be fixed. For sufficiently large $R,$ we choose $n:=\lfloor \log R\rfloor,$ with $\lfloor a \rfloor$ being the integral part of $a \in \mathbb{R},$ and 
 estimate the terms on the right-hand side of (15). Employing (10), we find that 
 $$|J_n(R)| \leq C e^{-(t-T)2n \pi}\int_{1/R}^R \frac{du}{u} = Ce^{-(t-T)2n\pi} \log R^2,$$
 and so by our choice of $n$ the term $J_n(R)$ approaches $0$ as $R \rightarrow \infty.$
 
 Next we estimate $I_k(\rho).$ Using again (10), we conclude that
 $$|I_k(\rho)| \leq C\int_0^{2 \pi} e^{-(t-T)(2k\pi + \theta)}d\theta = \frac{C}{t-T}e^{-2k\pi(t-T)}\bigg(1 - e^{-2\pi(t-T)}\bigg).$$
 Thus for an integer $n_0 \in (0, n-1)$ we obtain 
 \begin{eqnarray*}
 \sum_{k=n_0}^{n-1}|I_k(1/R) - I_k(R)| &\leq& \frac{2C}{t-T}\bigg(1- e^{-2\pi (t-T)}\bigg)\sum_{k=n_0}^\infty  e^{-2k\pi(t-T)} \\
 &=& \frac{2C}{t-T}e^{-2\pi n_0 (t-T)}.
 \end{eqnarray*}
 Now, the right-hand side can be made smaller than $\epsilon$ by fixing $n_0$ as a sufficiently large integer, which is possible since $R$ and so $n$ tends to infinity. Finally, by (11) we have
 $$\lim_{R \rightarrow \infty} \sum_{k=0}^{n_0-1}(I_k(1/R) - I_k(R)) = 0.$$
 Altogether, we have proved that for $t>T,$
 $$\lim_{R \rightarrow \infty}\left|\int_{1/R}^R g(u) e^{it\log u}\frac{d u}{u}\right| < \epsilon,$$
 for any $\epsilon >0.$ This implies that $[f]^\wedge_{M^2}(c+it) = 0.$
 When $t < -T,$ we choose $n = - \lfloor\log R\rfloor,$ sum in (15) from $n$ to $-1$ and proceed analogously. $\Box$
 
 \section{Bernstein type inequality with Mellin derivatives}
 
An important property of Mellin bandlimited functions is the following Bernstein inequality, involving Mellin derivatives
\begin{Theorem}
For any $f \in B^2_{c, T}$ one has $f \in C^\infty(\mathbb{R}^+)$ and moreover
$$\|\Theta_c^rf\|_{X^2_c} \leq (2T)^{r}\|f\|_{X^2_c}.$$
\end{Theorem}
\vskip0,3cm
A proof can be found in \cite{BJ3} for functions in $B^1_{c, T},$ but it also holds for functions in $B^2_{c, T},$ 
using formula (8.7) in \cite{BJ} and Lemma 2.7 in \cite{BJ4}. 
\vskip0,4cm
However, we give here a sharp version with $T^r$ in place of $(2T)^r$ directly for 
functions belonging to the space $\widetilde{B}^2_{c, T}.$ 
This leads to an alternative proof of the essential part of the Paley-Wiener theorem for Mellin transforms, namely Theorem 3, without employing 
Riemann surfaces. "Real" Paley-Wiener theorems for other kinds of integral transform can be found in \cite{AN}, \cite{AD}. 
\vskip0,3cm
We begin with the following
\begin{Lemma}
Let $f \in W^{2,r}_c,$ $r \in \mathbb{N}.$ Then, putting $h(x)= e^{cx}f(e^x), \quad x \in \mathbb{R},$ we have
\begin{eqnarray*}
\|h^{(r)}\|_2 = \|\Theta^r_cf\|_{X^2_c}
\end{eqnarray*}
\end{Lemma}
{\bf Proof.} We begin with $r=1.$ We have easily 
\begin{eqnarray*}
\|h'\|_2^2 &=& \int_{-\infty}^\infty |ce^{cx}f(e^x) + e^{(c+1)x}f'(e^x)|^2 dx \\
&=& \int_0^\infty |cu^c f(u) + u^{c+1}f'(u)|^2 \frac{du}{u}\\
&=& \int_0^\infty |u^c(cf(u) + uf'(u))|^2\frac{du}{u} \\&=&
\int_0^\infty |u^c \Theta_cf(u)|^2\frac{du}{u} = \|\Theta_cf\|^2_{X^2_c}.
\end{eqnarray*}
In the general case, we use the representation of the Mellin derivatives in terms of the Stirling numbers $S_c(r,k)$ of the second kind (see \cite{BJ},  Lemma 9)
$$(\Theta_c^rf)(u) = \sum_{k=0}^{r} S_c(r,k)u^k f^{(k)}(u) \quad \quad (u \in \mathbb{R}^+),$$
obtaining
\begin{eqnarray*}
\|\Theta^r_cf\|_{X^2_c}^2 &=& \int_0^\infty u^{2c}|\Theta^r_cf(u)|^2\frac{du}{u}\\
&=& \int_0^\infty u^{2c} \bigg|\sum_{k=0}^r S_c(r,k)u^k f^{(k)}(u) \bigg|^2\frac{du}{u},
\end{eqnarray*}
and so putting $u=e^x,$ $x \in \mathbb{R},$ the last integral is exactly $\|h^{(r)}\|^2_2.~~\Box$
\vskip0,4cm
\begin{Theorem}
[Bernstein inequality] Let $f \in \widetilde{B}^2_{c,T}.$ Then 
\begin{eqnarray}
\|\Theta_c^rf\|_{X^2_c} \leq T^r \|f\|_{X^2_c}.
\end{eqnarray}
\end{Theorem}
{\bf Proof.} We begin with the case $r=1.$ Since $f\in \widetilde{B}^2_{c,T},$ the function $h$ defined in Lemma 6 has an extension to an entire function of exponential type $T>0.$ 
One can assume $T=1,$ and deduce the general case by scaling the argument. Thus, using the proof of Proposition B' in \cite{BR} 
one has the representation
$$h'(0) = \frac{4}{\pi^2} \sum_{k \in \mathbb{Z}} (-1)^k \frac{h((2k+1)\pi/2)}{(2k+1)^2},$$
and the above equation holds for $h(t +\cdot)$ for any $t\in \mathbb{R},$ since the Bernstein spaces are shift-invariant. Hence
$$h'(t) = \frac{4}{\pi^2} \sum_{k \in \mathbb{Z}} (-1)^k \frac{h(t+(2k+1)\pi/2)}{(2k+1)^2}\quad\quad (t \in \mathbb{R}).$$
Applying the $L^2-$norm to both sides of this equation, and noting that the $L^2-$ norm is shift invariant, we get
$$\|h'\|_2 \leq \frac{4}{\pi^2}\sum_{k \in \mathbb{Z}}\frac{\|h\|_2}{(2k+1)^2} = \|h\|_2.$$
Replacing $h(s)$ with $h(s/T),$ we finally obtain 
$$\|h'\|_2 \leq T\|h\|_2.$$
Employing now Lemma 6, we obtain the assertion for $r=1.$ 
The general case follows by using an analogous argument, taking into account that the function $h^{(r)}$ is also of exponential type $T.$ Indeed, using the Cauchy integral theorem for the derivatives, one has, for $z \in \mathbb{C},$
$$h^{(r)}(z) = \frac{r!}{2 \pi i}\int_{{\cal C}} \frac{h(s)}{(s-z)^{r+1}}ds,$$
where ${\cal C}$ is the circle with center $z$ and radius 1. Then, since $h$ is of exponential type $T,$ we have for $z = x +iy,$
\begin{eqnarray*}
&&|h^{(r)}(z)| \leq \frac{r!}{2 \pi}\int_0^{2 \pi} |h(z +e^{it})|dt\\
&\leq & \frac{r! M}{2 \pi} e^{T|y|}\int_0^{2 \pi} e^{T|\sin t|}dt := H e^{T|y|}.
\end{eqnarray*}
Thus, we can proceed as before with $h^{(r)}$ in place of $h$ and using Lemma 6. $\Box$
\vskip0,4cm
By Theorem 5 we  now can give an alternative proof of the non-trivial part of the Paley-Wiener theorem, using real variable arguments. 
We follow a method similar to that given in \cite{AN} and the following definition of the support for the Mellin transform (see e.g. the general definition given in \cite{HIG}): for $f \in X^2_c,$ the support of $[f]^\wedge_{M^2},$ written $\mbox{Supp} [f]^\wedge_{M^2},$ is the set of all points $t \in \mathbb{R}$ that possess no neighbourhood throughout which $[f]^\wedge_{M^2} = 0,$ almost everywhere. The support is always a closed set. The next theorem shows that if $f \in \widetilde{B}^2_{c,T}$ then $\mbox{Supp} [f]^\wedge_{M^2}$ is contained in the interval $[-T,T].$
\vskip0,3cm
\noindent
\begin{Theorem}  $\widetilde{B}^2_{c,T} \subset B^2_{c,T}.$
\end{Theorem}
{\bf Proof}. By Theorem 5, we have that for every $r \in \mathbb{N},$ $\Theta^r_cf \in X^2_c.$ Moreover 
\begin{eqnarray}
\limsup_{r \rightarrow \infty} \|\Theta^r_cf\|_{X^2_c}^{1/r} \leq T 
\limsup_{r \rightarrow \infty}\|f\|_{X^2_c}^{1/r} = T.
\end{eqnarray}
Let now $t_0 \in \mbox{Supp}~[f]^\wedge_{M^2},$ with $t_0\neq 0,$ and let $\varepsilon >0$ be such that $|t_0| - \varepsilon >0.$ 
By Lemmas 4 and 5 we obtain
\begin{eqnarray*}
\|\Theta^r_cf\|^2_{X^2_c} &=& \frac{1}{2\pi}\|[\Theta^r_cf]^\wedge_{M^2}(c+i\cdot)\|^2_{L^2(\{c\}\times i \mathbb{R})}\\
&=& \frac{1}{2\pi}\|(-it)^r [f]^\wedge_{M^2}(c+i\cdot)\|^2_{L^2(\{c\}\times i \mathbb{R})}\\
&=& \frac{1}{2\pi}\int_{-\infty}^\infty |t|^{2r} |[f]^\wedge_{M^2}(c +i t)|^2 dt\\
&\geq& \frac{(|t_0| - \varepsilon)^{2r}}{2\pi}\int_{|t| > |t_0| -\varepsilon}|[f]^\wedge_{M^2}(c +i t)|^2 dt.
\end{eqnarray*}
Therefore, putting $H_\varepsilon:= \int_{|t| > |t_0| -\varepsilon}|[f]^\wedge_{M^2}(c +i t)|^2 dt,$ we have $H_\varepsilon \neq 0,$ and  
$$(|t_0| - \varepsilon)^{2r} \leq \frac{2 \pi}{H_\varepsilon}\|\Theta^r_cf\|^2_{X^2_c},$$
and then
$$|t_0| - \varepsilon \leq \bigg(\frac{2 \pi}{H_\varepsilon}\bigg)^{1/2r}\|\Theta^r_cf\|^{1/r}_{X^2_c}.$$
Finally by (17) passing to the lim sup we get $|t_0| \leq T.~~\Box$

\section{Some applications}

Two of the main applications of our Paley-Wiener theorem in the Mellin transform setting are the Exponential Sampling Formula and the Mellin Reproducing Kernel Formula. The first for functions in $B^{2}_{c, \pi T}$ was proved in \cite{BJ4}, while the second formula was given in \cite{BBM1}, where it was also proved that the two results are equivalent, in the sense that one result can be deduced from the other. Using now Theorem 3 we immediately can give versions of the above results directly for functions $f$ belonging to the Mellin-Bernstein space $\widetilde{B}^{2}_{c, T},$ with $T=\pi \sigma,$ $\sigma >0.$ 

Here we give the statements of the two results.

\begin{Theorem}[Exponential Sampling Formula]
If $f \in \widetilde{B}^2_{c, \pi \sigma}$ for some $c \in \mathbb{R},$ and $\sigma>0,$ then the series 
$$x^c \sum_{k=-\infty}^{\infty} f(e^{k/\sigma})\mbox{\rm lin}_{c/\sigma}(e^{-k}x^\sigma)$$
is uniformly convergent in $\mathbb{R}^+,$ and one has the representation
$$f(x) = \sum_{k=-\infty}^{\infty} f(e^{k/\sigma})\mbox{\rm lin}_{c/\sigma}(e^{-k}x^\sigma) ~~(x \in \mathbb{R}^+).$$
\end{Theorem}
\vskip0,4cm
The $\mbox{lin}_c-$function for $c\in I\!\!R,$ $\mbox{lin}_c : \mathbb{R}^+ \rightarrow \mathbb{R},$ is defined, for $x \in \mathbb{R}^+ \setminus \{1\},$ by
\begin{eqnarray} 
\mbox{lin}_c(x) = \frac{x^{-c}}{2 \pi i} \frac{x^{\pi i} - x^{-\pi i}}{\log x} = \frac{x^{-c}}{2 \pi }\int_{- \pi}^\pi x^{-it}dt,
\end{eqnarray}
with the continuous extension $\mbox{lin}_c(1) := 1,$ thus $\mbox{lin}_c(x) = x^{-c} \mbox{sinc} (\log x).$
\vskip0,3cm
It is clear that $\mbox{lin}_c \not \in X_{\overline{c}}$ for any $\overline{c}.$ However, it belongs to the space $X^2_c$ and its Mellin transform in $X_c^2-$sense is given by 
$$M^2[\mbox{lin}_c](c+it) = \chi_{[-\pi, \pi]}(t),~~t \in \mathbb{R}.$$ 
Here, $\chi_I$ denotes the characteristic function of the set $I.$
\vskip0,2cm
Note that integration of the formula of Theorem 7 over $\mathbb{R}^+$ gives a quadrature formula that was inverstigated in \cite{SCH}.

\begin{Theorem}[Mellin Reproducing Kernel Formula]
For $f \in \widetilde{B}^2_{c,\pi \sigma},$ with $c \in \mathbb{R}$ and $\sigma>0,$ we have 
$$f(x) = \sigma\int_0^{\infty} f(y)\mbox{\rm lin}_{c/\sigma}\bigg((\frac{x}{y})^\sigma\bigg)\frac{dy}{y},$$
 the integral being absolutely convergent in $\mathbb{R}^+.$
\end{Theorem}

\end{document}